\documentclass[11pt]{article}

\usepackage{amsmath}
\usepackage{amssymb}
\usepackage{amsfonts}
\usepackage[utf8]{inputenc}
\usepackage[english]{babel}
\usepackage{mathrsfs}
\usepackage{algorithm,algorithmicx}
\usepackage{xspace}
\usepackage[margin=2.5cm]{geometry}
\usepackage{bm}

\usepackage{hyperref}
\hypersetup{
    unicode = false,
    pdftoolbar = true,
    pdfmenubar = true,
    pdffitwindow = true,
    pdftitle = {hierarchy},
    pdfauthor = {Alarie et al.},
    pdfsubject = {hierarchy},
    pdfkeywords = {hierarchy},
    pdfnewwindow = true,
    colorlinks = true,
    linkcolor = blue,
    citecolor = blue,
    filecolor = black,
    urlcolor = blue,
    breaklinks = true
}

\usepackage{cleveref}
\usepackage{pbox}
\usepackage{tikz}
\usetikzlibrary{arrows,backgrounds,decorations,positioning,intersections,calc}

\def\R{{\mathbb{R}}}
\def\N{{\mathbb{N}}}

\def\boite{\phantom{.}\hfill{$\Box$}}
\newcommand{\mads}{{\sf Mads}\xspace}
\newcommand{\eb}{{\sf EB}\xspace}
\newcommand{\pb}{{\sf PB}\xspace}
\newcommand{\intr}{{\sf Int}\xspace}
\newcommand{\hier}{{\sf Hier}\xspace}
\newcommand{\nomad}{{\sf NOMAD}\xspace}

\makeatletter
\renewcommand{\ALG@name}{Algo.}
\makeatother

\usepackage{algpseudocode}

\newtheorem{Def}{Definition}[section]

\newtheorem{Prop}[Def]{Proposition}

\newcommand{\Demo}{\noindent {\bf Proof.}~}

\usepackage{color}
\definecolor{Red}{rgb}{1,0,0}
\definecolor{Blue}{rgb}{0,0,1}

\newcommand{\rinf}{\overline{\mathbb{R}}}

\newtheorem{algolocal}{Algorithm}[section]
\newenvironment{algo}[1]
	{	\begin{algolocal} {\em \bf #1} \hrulefill \end{algolocal}
		\begin{list}{}{} \item \vspace*{-2mm}			}
	{	\end{list} \vspace*{-3mm}\noindent \hrulefill \\			}

\title{
Hierarchically constrained blackbox optimization
\thanks{This work has been funded by a MITACS Elevate grant in collaboration with Hydro-Québec.}}
\author{
\href{mailto:alarie.stephane@hydroquebec.com}{St\'ephane Alarie}
\thanks{
    Hydro-Qu\'ebec and IREQ,
    \href{mailto:alarie.stephane@hydroquebec.com}{\tt alarie.stephane@hydroquebec.com}.
}
\and
\href{mailto:Charles.Audet@gerad.ca}{Charles Audet}\thanks{
    {GERAD} and Polytechnique Montr\'eal,
          \href{https://www.gerad.ca/Charles.Audet}{\tt www.gerad.ca/Charles.Audet}.
}
\and
\href{mailto:jacquot.paulin@gmail.com}{Paulin Jacquot}\thanks{
EDF,
\href{mailto:jacquot.paulin@gmail.com}{\tt jacquot.paulin@gmail.com}.
}
\and
\href{mailto:Sebastien.Le.Digabel@gerad.ca}{S\'ebastien {Le~Digabel}}\thanks{
{GERAD} and Polytechnique Montr\'eal,
\href{https://www.gerad.ca/Sebastien.Le.Digabel}{\tt www.gerad.ca/Sebastien.Le.Digabel}.
} 
}
\date{\today}

\begin{document}

\maketitle

\newcommand{\pj}[1]{{\color{purple}#1}}

\begin{abstract}
In blackbox optimization,
    evaluation of the objective and constraint functions is time consuming.
In some situations, constraint values may be evaluated independently or sequentially.
The present work proposes and compares two strategies to 
    define a hierarchical ordering of the constraints and to
    interrupt the evaluation process at a trial point 
    when it is detected that it will not improve the current best solution.
Numerical experiments are performed on a closed-form test problem. \\

\noindent
{\bf Keywords:}
Blackbox Optimization, Derivative-Free Optimization, Constrained Optimization.
\end{abstract}


\section{Introduction}

Consider the constrained optimization problem
$$(P) \qquad \qquad    \begin{array}{rl}\displaystyle 
    \min_{x\in {\cal X}} & f(x) \\
      \mbox{s.t.} & c_1(x) \leq 0 \\
        & c_2(x) \leq 0 \\
        &   \qquad \vdots  \\
        & c_m(x) \leq 0
    \end{array}$$
in which $\cal X$ represents the domain of the objective $f: {\cal X} \rightarrow \rinf$ and of the constraint functions  $c_j:{\cal X} \rightarrow \rinf$, for $j \in \{1,2,\ldots, m\}$.
The domain $\cal X$ is usually $\R^n$, $\N^n$ or some bound-constrained subset of one of these spaces.
The entire feasible region is compactly denoted by
$$ \Omega \ = \ \left\{ ~ x \in {\cal X} ~ | ~ c_j(x) \leq 0 \ \forall j \in \{1,2,\ldots, m\} ~\right\}.$$

The present work considers the situation where the functions defining $(P)$ 
are provided as \emph{blackboxes}~\cite{AuHa2017}, usually in the form of computer codes,
 but more specifically, cases in which each constraint $c_j$ and the objective $f$
    may be evaluated independently or sequentially.
We adopt the convention that $c_j(x)=+\infty$ if the evaluation of $c_j(x)$ could not be terminated, for example if the simulation process crashed.
A similar convention holds for the objective function $f$.

More precisely, we consider the case where the $j$-th blackbox returns 
    the value of the $j$-th constraint $c_j$ for $j$ ranging from $1$ to $m$,
    and the $(m+1)$-th blackbox returns the objective function $f$.
When studying a trial point $x \in \R^n$,
    a decision may then be taken after each constraint evaluation
    to decide whether or not to spend more effort in 
    evaluating the remaining functions.

The objective of this work is to solve Problem~$(P)$
    while reducing as much as possible the overall computational effort.
The motivating problem is the PRIAD project at Hydro-Québec  TransÉnergie~\cite{PRIAD_KoMeCoGaVoAlDeBl2021}
    to optimize the maintenance strategies of its electrical power grid.
The optimization model is still being developed,
    but each evaluation will go through a sequence of five blackboxes given the frequencies of preventive maintenance to be performed on grid equipment:
    i- feasibility with the workload capacity;
    ii- resulting equipment failure probabilities;
    iii- scenario generation of unavailable equipment;
    iv- power flow simulations to quantify energy delivery interruptions;
    v- risk and cost assessment of load shedding and corrective maintenance. 
The first four are related to constraints and the last to the objective function.
They are chained so that the outputs of one are inputs of the next.
Computational time also increases significantly from blackbox to blackbox.
While the first two will take a few seconds to complete, the fourth will require days.

We propose two different two-phase procedures, in
 which the first phase is dedicated to finding a feasible solution,
 and the second phase is devoted to the optimization of $(P)$.
Both procedures exploit the flexibility of the Mesh Adaptive Direct Search (\mads) 
    algorithm~\cite{AuDe2006} designed for constrained blackbox optimization.

This document is composed of two main sections.
Section~\ref{sec-algos}
    describes a pair of algorithmic procedures to interrupt evaluations
and
Section~\ref{sec-exp}
    illustrates the procedures on an analytical 
    tension/compression spring design problem.
Concluding remarks and potential future work close the paper in Section~\ref{sec-Discussion}.

\section{Interruption of futile blackbox evaluations}
\label{sec-algos}
Subsections~\ref{sec-interruptible} and~\ref{sec-hierarchical} describe 
algorithmic procedures to interrupt the sequence of constraint evaluations.
Both procedures rely on the \mads framework described in the next subsection.

\subsection[The Mesh Adaptive Direct Search algorithm (MADS)]
           {The Mesh Adaptive Direct Search algorithm (\mads)}
\mads is a class of direct search algorithms designed to solve Problem~$(P)$.
It was first proposed in the paper~\cite{AuDe2006},
    but has undergone many improvements.
The latest description is found in the recent paper~\cite{AuLeDTr2018}.
The simplest form of this class of algorithms is presented in 
the textbook~\cite{AuHa2017}.

\mads is an iterative algorithm that attempts at each iteration 
    to generate a trial point that would be an improvement over 
    the current best-known solution called the incumbent.
\mads is initiated with one trial point $x^0 \in {\cal X}$
    that is not required to satisfy the constraints $c_j$.
Each trial point $x \in \R^n$ is sent to the blackbox for the evaluation
    of the constraints $c_j(x)$ and objective $f(x)$ function values.
The incumbent solution at iteration $k$ is denoted $x^k$.

Two mechanisms to handle the constraints are proposed in the blackbox optimization literature.
Both use the {\em constraint violation function}
$$ h(x)  \ = \ 
    \left\{ \begin{array}{cl} \displaystyle 
        \sum_{j=1}^m \max(c_j(x),0)^2 & \mbox{ if } x \in {\cal X}\\
        +\infty & \mbox{ otherwise}
        \end{array}\right.$$
inspired from filter methods~\cite{FlLe02a} and adapted in~\cite{AuDe09a}.
The constraint violation function value $h(x)$ is zero if and only if $x$ 
 belongs to the feasible region $\Omega$.

The first mechanism, referred to as the Extreme Barrier (\eb) approach, is a two-phase method solving a pair 
    of optimization problems without general constraints.  
The first phase attempts to find a feasible solution
    by minimizing the constraint violation function $h$ from the starting point $x^0$.
It is only invoked if $x^0$ is infeasible with respect to the constraint $c(x) \leq 0$.
At the starting point, the constraint violation function value $h(x^0)$ is finite since $x^0 \in {\cal X}$.
If the first phase fails to generate a feasible point, 
    then Problem~$(P)$ is declared infeasible and the process terminates.
Otherwise, as soon as a feasible point $\bar x$ is generated, 
 the first phase is stopped.
The second phase is then launched, and consists in minimizing the extreme barrier function
$$\displaystyle f_\Omega(x) := \left\{ \begin{array}{cll}
    f(x) &~\quad~ &\mbox{if}~x \in \Omega, \\
    +\infty && \mbox{if}~x \notin \Omega
    \end{array}\right.$$
from the feasible starting point $\bar x$.

The second mechanism to handle the constraints is called the Progressive Barrier (\pb)~\cite{AuDe09a},
    and consists in analyzing the bi-objective problem of minimizing 
    both the objective function $f$ and the constraint violation function $h$.
In order to apply \pb,
  all the constraint functions $c_j$ need to be evaluated at the trial points.
We next introduce two novel procedures that allow incomplete evaluations of the constraint functions.

\subsection{The interruptible extreme barrier procedure}
\label{sec-interruptible}
A first approach to solve Problem~$(P)$ 
    consists in launching \mads with the two-phase extreme barrier,
    but to simply interrupt the sequence of constraint evaluations as soon as one is violated.

The first phase of the approach,
 the feasibility phase, is again to minimize the constraint violation function $h$ 
 as for \eb above, from a user-provided starting point $x^0$,
    but to terminate the evaluations of the constraints as soon as the cumulative constraint violation function exceeds the incumbent value.
Indeed, the partial sums  $$s_\ell(x) := \sum_{j=1}^\ell \max(c_j(x),0)^2$$ appearing in $h$ 
    are monotone increasing with respect to the index $\ell$ 
    and therefore, the incumbent solution $x^k$ is replaced by a trial point $t \in {\cal X}$ if and only if $$s_\ell(t) < h(x^k) 
    \quad \mbox{ for every index } \quad \ell \in \{ 1,2,\ldots,m\}.$$
   
The second phase solves Problem~$(P)$
    from the feasible starting point $\bar x$ produced by the first phase (if it succeeded).
Again, the evaluations of the constraints is interrupted as soon as a violation is detected.
Algorithm~\ref{Algo-2 phase method}
    gives the pseudo-code of the two-phase
    interruptible extreme barrier procedure,
    denoted (\intr).

\begin{algo}{
	\label{Algo-2 phase method}
	\sf 	Two-phase interruptible extreme barrier (\intr)}
Given the objective function $f : \R^n \rightarrow \rinf$,
    constraints $c_j : \R^n \rightarrow \rinf$
    for $j \in \{1,2,\ldots,m\}$ and starting point $x^0 \in {\cal X}$.\\
{\sf 1. Feasibility phase}\\
	\hspace*{5mm} \begin{tabular}[t]{|l}
	   Apply \mads with the extreme barrier 
	    from  $x^0 \in {\cal X}$
	to solve $\displaystyle  \min_{x}  h(x)$\\
	\hspace*{5mm} \begin{tabular}[t]{|l}
	For a trial point $t \in {\cal X}$,
	    evaluate $\{c_\ell(t)\}_{\ell =1}^m$ in sequence,\\
\hspace*{5mm}	    but reject $t$ and interrupt evaluation as soon as $s_\ell(x) \geq h(x^k)$
	    for some index $\ell$\\
	    	Interrupt \mads and terminate the feasibility phase as soon as a \\
\hspace*{5mm} feasible point $\bar x \in \Omega$ is generated and go to Phase {\sf 2}\\

    \end{tabular}	 \\   
	If \mads failed to generate a feasible point, \\
\hspace*{5mm} 	terminate by
	concluding that no point in $\Omega$ was found
\end{tabular} \\ \\
{\sf 2. Optimization phase}\\
\hspace*{5mm} \begin{tabular}[t]{|l}
	   Apply \mads with the extreme barrier (\eb)
	    starting from $\bar x \in \Omega$
	to solve Problem~$(P)$\\
	\hspace*{5mm} \begin{tabular}[t]{|l}
	For a trial point $t \in {\cal X}$,
	    evaluate $\{c_\ell(t)\}_{\ell =1}^m$ in sequence,\\
\hspace*{5mm}	    but reject $t$ and interrupt  evaluation as soon as $c_\ell(t) > 0$
for some index $\ell$
\end{tabular}\\
\end{tabular} \\
\end{algo}


The next result gives sufficient conditions to ensure that \mads with \eb and \intr behave identically (in the sense that they generate the same sequence of iterates),
    however, the latter requires a lower computational effort.
The main condition is that no dynamic surrogate model is used by the algorithm,
    because the models are constructed by using previously evaluated function values.
Therefore, more points will be available to construct models to the \eb procedure than to the \intr procedure.

\begin{Prop}
\label{Prop:same}
\mads with \eb and \intr produce the same sequence of iterates 
    when launched on the same Problem~$(P)$ from the same starting point
    and without using dynamic surrogate models.
\end{Prop}
\Demo
If \eb and \intr start an iteration from the same incumbent solution and with the same algorithmic parameters and random seeds, both will create the same list of trial points when no dynamic models are used.

First, consider feasibility phase.
If the starting point $x^0$ satisfies all the constraints,
    then both the \intr and \eb approaches conclude the phase with $\bar x = x^0$.
Otherwise,
    both approaches attempt to minimize $h(x)$ on ${\cal X}$.
At iteration $k$, 
    any trial point $t$ outside of ${\cal X}$ is rejected and 
    $x^{k+1}$ is set to $x^k$, the current incumbent solution.
If $t$ belongs to ${\cal X}$, 
    then there are two possibilities:
\begin{itemize}
\item[i-] If $h(t) < h(x^k)$,
        then both approaches will set $x^{k+1}$ to $t$;
\item[ii-] If $h(t) \geq h(x^k)$,
            then both approaches will set $x^{k+1}$ to $x^k$.
\end{itemize}
It follows that the feasibility phase of both approaches behave identically.

Second, consider the optimization phase.
Any trial point $t$ outside of $\Omega$ is rejected by both approaches and $x^{k+1}$ is set to $x^k$.
If $t$ belongs to $\Omega$, 
    then there are again two possibilities:
\begin{itemize}
\item[i-] If $f(t) < f(x^k)$,
        then both approaches will set $x^{k+1}$ to $t$;
\item[ii-] If $f(t) \geq f(x^k)$,
        then both approaches will set $x^{k+1}$ to $x^k$.
\end{itemize}
Once more, it follows that both optimization phases behave identically.
\boite

\subsection{The  hierarchical satisfiability extreme barrier procedure}
\label{sec-hierarchical}


A second approach to solve Problem~$(P)$ while reducing the computational burden 
    sequentially solves a collection of problems parameterized
    by the constraint index $j \in \{1,2,\ldots,m\}$,
    that minimize the function $c_j$ subject to the constraint with lesser indices
%
$$    (P_j) \qquad\qquad
    \begin{array}{rll}\displaystyle 
    \min_{x\in {\cal X}} & c_j (x) \\
      \mbox{s.t.} & c_1(x) \leq 0 \\
        & c_2(x)    \leq 0 \\
        &   \qquad \vdots  \\
        & c_{j-1}(x)  \leq 0
    \end{array}
$$

The feasibility phase starts by solving the unconstrained Problem~$(P_0)$
    from the user-provided starting point $x^0 \in {\cal X}$.
The algorithm then adds one constraint at a time,
 treats it with the extreme barrier and minimizes the next one
 by sequentially solving $(P_1)$ to $(P_m)$.
Evaluations are interrupted as soon as a constraint is violated, and
 the entire process is terminated when \mads is unable to
  find a feasible solution to one of  Problems~$(P_j)$.
The optimization phase is identical to that of the \intr approach.

The pseudo-code of the hierarchical  satisfyability  procedure called (\hier) is presented next.

\begin{algo}{
	\label{Algo-Hier}
	\sf 	Hierarchical satisfiability with the extreme barrier (\hier)}
Given the objective function $f : \R^n \rightarrow \rinf$,
    constraints $c_j : \R^n \rightarrow \rinf$
    for $j \in \{1,2,\ldots,m\}$ and starting point $x^0 \in {\cal X}$.\\
{\sf 1. Feasibility phase}\\
	\hspace*{5mm} \begin{tabular}[t]{|l}
	   For each index $j$ varying from $1$ to $m$, \\
	   \hspace*{5mm} \begin{tabular}[t]{|l}
	   Apply \mads with the extreme barrier (\eb)
	    from $x^{j-1} \in {\cal X}$
	to solve Problem~$(P_j)$ \\
	\hspace*{5mm} \begin{tabular}[t]{|l}
    Terminate \mads as soon as a feasible point is found and denote it $x^j$\\
    For a trial point $t$,
	    evaluate $\{c_\ell(t)\}_{\ell =1}^{j-1}$ in sequence,\\
\hspace*{5mm}	    but interrupt the evaluation process as soon as $c_\ell(t) > 0$ 
	    for some index $\ell$
    \end{tabular}\\
    Terminate the entire algorithm if no feasible point is found by \mads
    \end{tabular}
\end{tabular} \\ \\
{\sf 2. Optimization phase}\\
\hspace*{5mm} \begin{tabular}[t]{|l}
	   Apply \mads with \eb 
	    from the starting point  $x^m \in \Omega$
	to solve Problem~$(P)$\\
	\hspace*{5mm} \begin{tabular}[t]{|l}
	For a trial point $t$,
	    evaluate $\{c_\ell(t)\}_{\ell =1}^m$ in sequence,\\
\hspace*{5mm}	    but interrupt the evaluation process as soon as $c_\ell(t) > 0$ 
	    for some index $\ell$
\end{tabular}\\\end{tabular} \\
\end{algo}

The \hier approach does not produce the same sequence of trial points 
    than \mads with \eb nor \intr.
\section{Illustration on a Tension/Compression Spring Design problem}
\label{sec-exp}
We compare the interruptible and hierarchical approaches described above 
    with the baseline given by direct calls to the \nomad implementation~\cite{Le09b} of the \mads algorithm,
    using the extreme or progressive barrier to handle constraints. 
This results in four optimization algorithmic procedures:
\begin{itemize}
\item 
\intr:  the  interruptible extreme barrier approach
    from Section~\ref{sec-interruptible};
\item 
\hier: the  hierarchical satisfiability extreme barrier approach
    from Section~\ref{sec-hierarchical};
\item 
\eb: \mads with the extreme barrier to handle constraints~\cite{AuDe2006};
\item 
\pb : \mads with the progressive barrier to handle constraints~\cite{AuDe09a};
\end{itemize}



 The Tension/Compression Spring Design problem (TCSD) consists of minimizing the
 weight of a spring under mechanical constraints. 
 The problem was originally introduced in~\cite{belegundu1985study} and has  also been considered  in~\cite{garg2014solving} along with other nonlinear engineering problems.
Although it is not a blackbox optimization problem because closed-form expressions of all functions are known,
    we use it to illustrate the behaviour of the algorithmic procedures.

This problem has three bound-constrained variables and four nonlinear constraints. The design variables
define the geometry of the spring. The constraints concern limits on outside
diameter~\eqref{cs:outdiam}, surge frequency~\eqref{cs:frequency},
minimum deflection~\eqref{cs:deflection} and shear stress~\eqref{cs:sheerstress}:
\begin{subequations}
 \begin{eqnarray}
 \label{pb:TCSD}
   f(x) & =&    x_{1}^2 \,  x_{2} \, (x_{3}+2) \\
   c_1(x) &=& \frac1{1.5}\, (x_{1}+x_{2}) -1 \ \leq \ 0 \label{cs:outdiam}\\
    c_2(x) &=& 1- 140.45 \, \frac{x_{1}}{x_{2}^2 \,x_{3}} \ \leq 0\   \label{cs:frequency}\\   
    c_3(x) &=& 1- \frac{x_{2}^3 \, x_{3}}{71785 \, x_{1}^4} \ \leq \ 0 \label{cs:deflection} \\
    c_4(x) &=& \frac{4 \, x_{2}^2 - x_{1}\,x_{2}}{12566 (\, x_{1}^3\,x_{2} - x_{1}^4)} + \frac1{5108\,x_{1}^2} -1 \ \leq \ 0 \label{cs:sheerstress} 
\end{eqnarray}
\end{subequations}
This problem is difficult in two ways.
It is hard to find feasible solutions and it possesses a large number of local optima.

The best known solution, denoted $x^*$ hereafter, and the bounds on the variables,
denoted by $\ell$ and $u$, are given in Table~\ref{tab-spring}.
The objective function value of $x^*$ is $f(x^*)= 0.0126652$.
\begin{table}[htb!]
  \centering
  \begin{tabular}{|ll||rr|r|}
    \hline
    \multicolumn{2}{|l||}{variable} & $\ell$ &  $u$ & $x^*$  \\\hline
    $x_1$ & mean coil diameter &0.05 & 2.0 & 0.051686 \\
    $x_2$ & wire diameter & 0.25 &1.3 & 0.35666 \\
    $x_3$ & active coils length & 2.0 & 15.0 & 11.29231 \\\hline
  \end{tabular}
\caption{Bounds and best know feasible point for the TCSD problem.}
\label{tab-spring}
\end{table}

For illustrative purposes, 
    the \textit{costs} $(\tau_k)$
    of evaluating the functions is 
    related to the evaluation time and is set to be the \textit{number of multiplications and divisions} for each constraint.
Table~\ref{tab-costs-spring}
    summarizes these costs.
The constraints are ordered so that their costs increase with the index number.
\begin{table}[htb!]
\begin{center}
    \begin{tabular}[c]{|c||ccccc|}
     \hline
function &      $c_1$ & $c_2$& $c_3$& $c_4$ & $f$ \\\hline
cost $\tau_k$ &   1 & 4 & 8 & 14 &3 \\\hline
  \end{tabular}
\end{center}
\caption{Costs of evaluating the functions defining the TCSD problem.}
\label{tab-costs-spring}
\end{table}

In the numerical experiments below,
the stopping criteria consists of an overall budget of  10,000 multiplications and divisions required by the evaluation of the functions.

The four algorithmic procedures are applied to $40$ instances of the TCSD problem
    by randomly generating infeasible starting points within the bound-constrained domain
    using an uniform distribution.

Two series of tests are conducted.
The first one orders the constraints so that the least expensive are evaluated first
 and the second one orders the constraints by considering first the ones that are most likely to be violated.

\subsection*{Constraints ordered by increasing evaluation cost}

The four constraints are ordered as in Table~\ref{tab-costs-spring} and
 results are summarized in \Cref{tab:results_TCSD}.
The left part of the table shows the average cost to reach the first feasible trial point encountered by each procedures.
The columns $c_1$ to $c_4$ and $f$ indicate the average number of times that each function
    is called to reach the final objective function value.
Even though there are no explicit mechanism to interrupt the sequence
    of constraint evaluations with the algorithmic procedures \eb and \pb, 
    the number of evaluations is not constant.
For both procedures the table reveals that the objective $f$ was not evaluated as often as the constraints.
The reason is that there were trial points for which the constraint $c_4$ generated a division by zero
    and this lead to an immediate termination of the simulation without evaluating 
    the objective function $f$. 

\begin{table}[htb!]
  \centering
\begin{tabular}{|c||r||ccccc||cc|}
\cline{2-9}
\multicolumn{1}{c||}{~}&
\multicolumn{1}{c||}{Feasible}& 
\multicolumn{5}{c||}{Final number of evaluations} &
\multicolumn{2}{c|}{Final objective $f(x^*)$}\\
\hline
Proc. & Avg. cost  & $c_1$ & $c_2$& $c_3$& $c_4$ & $f$  & instances avg. & best instance\\  \hline \hline
\eb   & 1479.9 & 334.0 & 334.0 & 334.0 & 334.0 & 333.7 & 0.0129921 & 0.0126654 \\
\pb   & 2199.4 & 334.0 & 334.0 & 334.0 & 334.0 & 333.8 & 0.0131401 & 0.0126656 \\
\intr &  939.0 & 457.3 & 449.4 & 440.2 & 264.0 & 179.1 & 0.0129580 & 0.0126659 \\
\hier &  833.8 & 447.8 & 437.5 & 429.6 & 278.5 & 160.2 & 0.0133600 & 0.0126654 \\
\hline
\end{tabular}
\caption{Average costs for reaching feasibility and final objective function values within an allocated budget of 10,000 multiplications and divisions
with constraints ordered by evaluation cost.}
\label{tab:results_TCSD}
\end{table}

The right part of \Cref{tab:results_TCSD} shows the final objective function value after that the budget of 10,000 multiplications and division is spent.
One column gives the average value and the other gives the lowest objective function value over the $40$ instances.

With the \intr and \hier procedures, 
    the number of times that the functions are called decreases monotonically from $c_1$ to $f$.
The table shows that the number evaluations does not vary importantly for the constraints $c_1, c_2$ and $c_3$.
This is explained by the fact that $c_1$ and $c_2$ are strictly satisfied
    at the optimal solution $x^*$
    and consequently the evaluations are not often interrupted for most trial points.
The constraints  $c_3$ and $c_4$ are binding at $x^*$
    and thus a significant number of evaluations of $c_4$ and $f$ are spared.
The table also shows that the different procedures require very different costs to generate a first feasible point. 
On average, \hier finds a first feasible point almost 60\%
 faster than \pb.

\begin{figure}[htb!]
  \centering
  \includegraphics[width=0.8\textwidth]{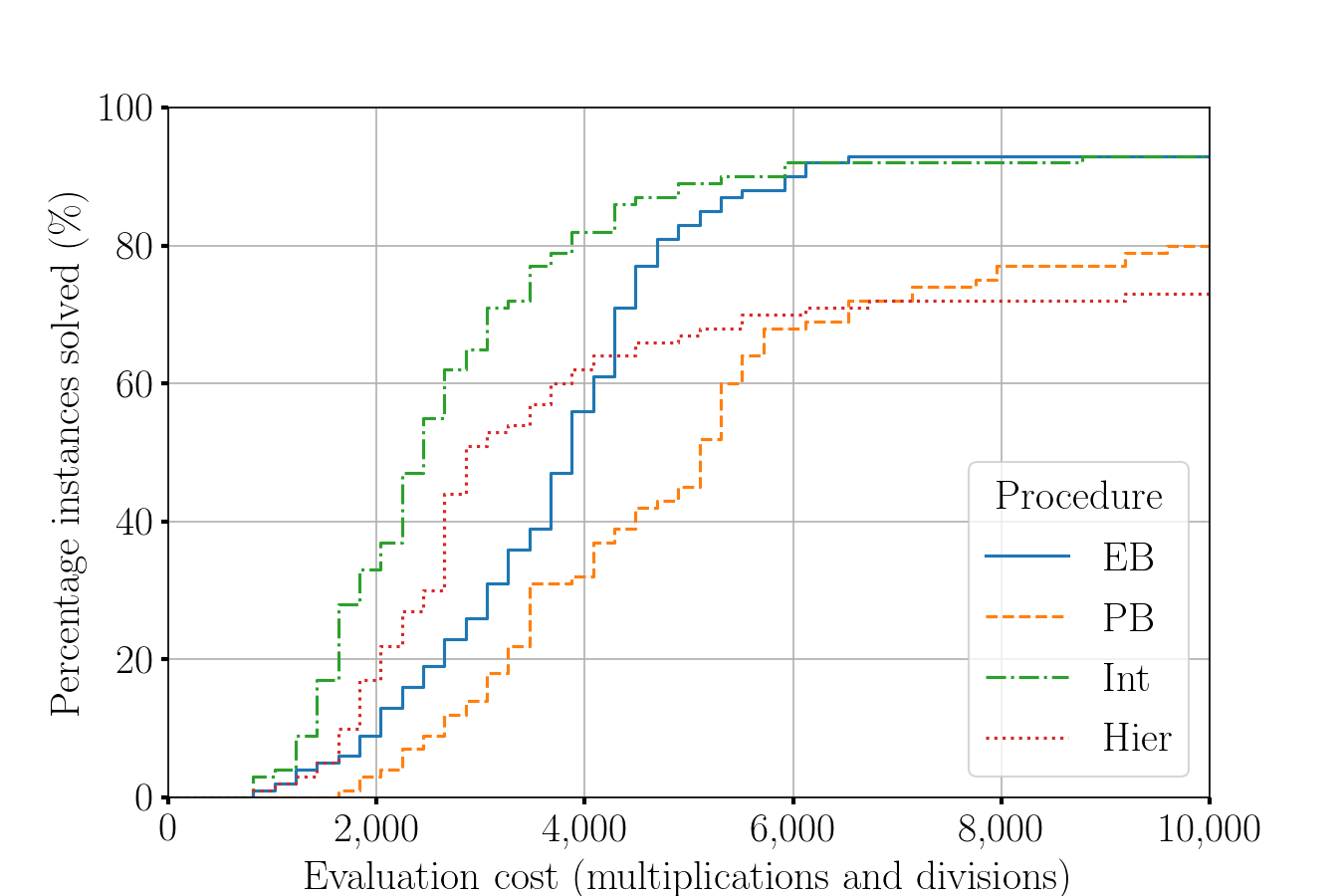}
  \caption{Data profiles showing performance of the different procedures on the TCSD problem, considering constraints in \textit{increasing} order of evaluation cost (criterion:  optimality gap of 5\%).}
  \label{fig:evo_proc_TCSD}
\end{figure}

Figure~\ref{fig:evo_proc_TCSD}
    compares the four algorithmic procedures
    with data profiles~\cite{MoWi2009} that show the proportion of instances solved 
    with a $5\%$ tolerance (i.e., a 5\% gap to the best known value $f(x^*)$) versus the evaluation cost.
None of the curves are connected to the ordinate axis because of the infeasible starting points.
Each curve starts when the procedure 
 has solved one of the 40 instances within $5\%$.
For \intr, \hier and \eb algorithmic procedures, the first instance is solved with the same effort, while
the \pb approach requires more multiplications and divisions.
This is coherent with the fact that \pb places more effort around infeasible points 
    with promising objective function values~\cite{AuDe09a}.
    
Over most of the duration of the optimization,
    the curve corresponding to the \intr procedure 
    dominates the others.
Close inspection of the figure
    reveals that at an evaluation cost of approximately 7,000,
    the \eb
    procedure slightly outperforms the \intr procedure.
This does not contradict Proposition~\ref{Prop:same} 
    because the construction and utilisation of dynamic quadratic models
    of the objective and constraint functions was enabled, as this is the default 
    option in \nomad.
This behaviour suggests that the use of quadratic models on this simple TCSD problem helps to accelerate the convergence when the basin containing a local optima is found.
When the number of iterations is low, the models are not as useful, and both the \intr and \hier procedures benefit from interrupting the calls to the simulation.

\subsection*{Constraints ordered by  empirical infeasibility}

There are situations where one knows which constraints are difficult to satisfy,
     which ones are easily satisfied and which ones are not.
For example,
    Rio Tinto's engineers
    know that maintaining a small flooding risk for the Kemano
    hydroelectric system is much more challenging than 
    ensuring continuous smelter operations~\cite{CoPa2019}.
For the TCSD problem, 
    we simulated that knowledge by evaluating the functions at 5,000 points using Latin Hypercube sampling~\cite{McCoBe79a} in the domains delimited by the bounds listed in Table~\ref{tab-spring}.
Table~\ref{Tab:feasible} gives the proportion of sample points that are feasible for each constraint.

\begin{table}[htb!]
\begin{center}
    \begin{tabular}[c]{|c||cccc|}
     \hline
Constraint    &      $c_3$ & $c_1$& $c_2$& $c_4$  \\\hline
Feasibility (\%)    &   2 & 34 & 99 & 99 \\\hline
  \end{tabular}
\end{center}
\caption{Proportion of feasible points for each constraint when sampling 5,000 points, ranked in increasing order.}
\label{Tab:feasible}
\end{table}

Ordering the constraints so that the most likely to be infeasible are evaluated first 
    will potentially trigger the interruption at a lower computational cost.
Results of the performances of the different procedures on the TCSD problem are shown in Table~\ref{tab:results_TCSD_DECREASING_order_constraints}.

\begin{table}[htb!]
  \centering
\begin{tabular}{|c||r||ccccc||cc|}
\cline{2-9}
\multicolumn{1}{c||}{~}&
\multicolumn{1}{c||}{Feasible}& 
\multicolumn{5}{c||}{Final number of evaluations} &
\multicolumn{2}{c|}{Final objective $f(x^*)$}\\
\hline
Proc. & Avg. cost  & $c_3$ & $c_1$& $c_2$& $c_4$ & $f$  & instances avg. & best instance\\  \hline \hline
\hline
\eb & 1264.8 & 334.0 & 334.0 & 334.0 & 334.0 & 333.8 & 0.0130931 & 0.0126653\\
\pb & 1934.4 & 334.0 & 334.0 & 334.0 & 334.0 & 333.7 & 0.0131127 & 0.0126656 \\
\intr & 863.8 & 485.4 & 290.7 & 290.3 & 280.2 & 195.2 & 0.0130218 & 0.012668 \\
\hier & 863.1 & 473.1 & 309.9 & 309.9 & 296.5 & 165.4 & 0.0133596 & 0.0126653
 \\\hline
\end{tabular}
\caption{Average costs for reaching feasibility and final objective function values within an allocated budget 10,000 multiplications and divisions, 
with constraints ordered by empirical infeasibility.}
\label{tab:results_TCSD_DECREASING_order_constraints}
\end{table}

\begin{figure}[htb!]
  \centering
  \includegraphics[width=0.8\textwidth]{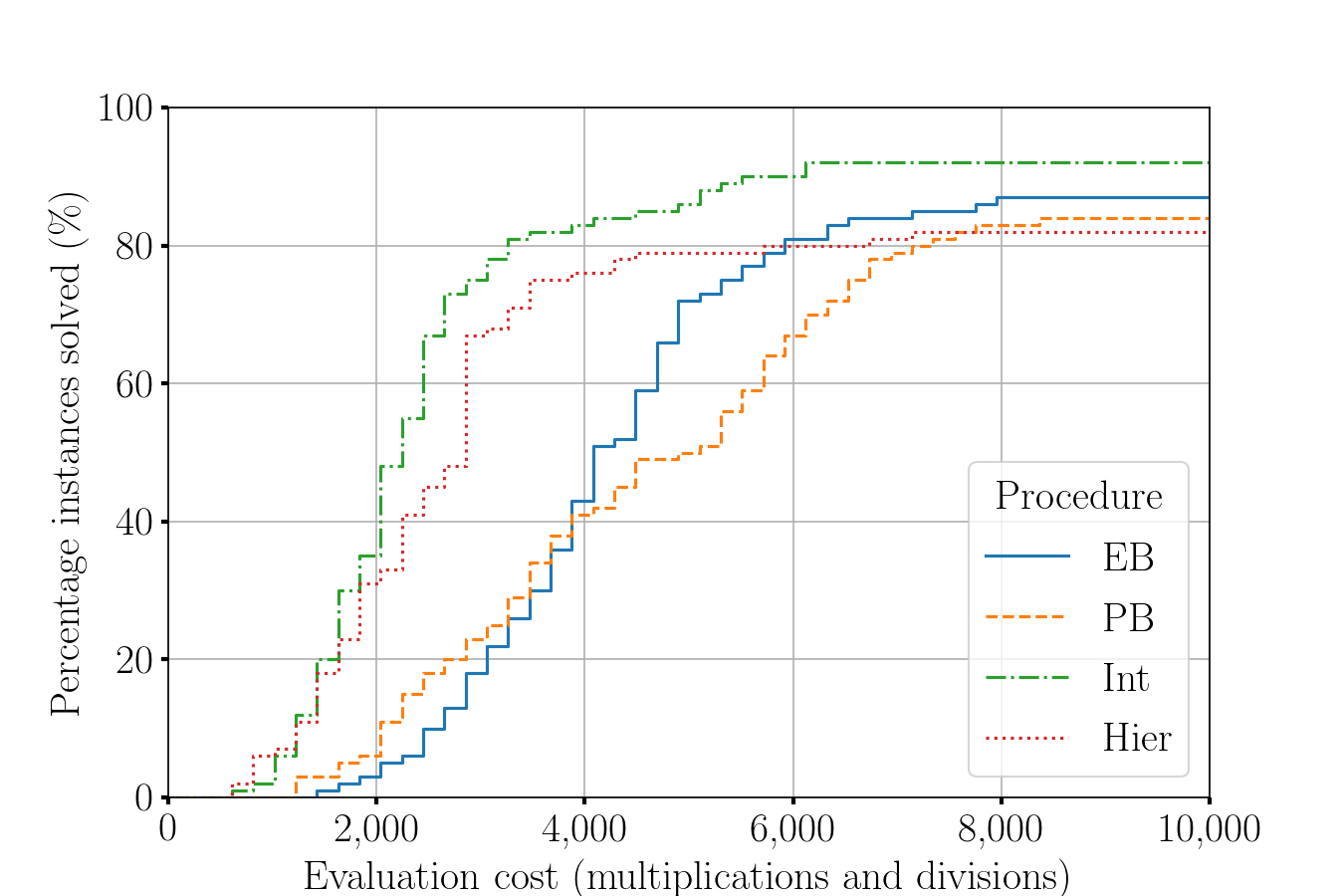}
  \caption{Data profiles showing performance of the different methods on the TCSD problem, with constraints ordered by empirical infeasibility (criterion:  optimality gap of 5\%).}
  \label{fig:evo_proc_TCSD_DECREASING_order_constraints}
\end{figure}

A first observation is that results of the \eb and \pb procedures differ from those listed in Table~\ref{tab:results_TCSD}.
This discrepancy is due to the random selection of starting points.
A second observation is that the average cost for the \intr procedure drops 
 by approximately 8\% for the generation of a first feasible solution,
 and increases by 3.5\% for \hier.
The average and best final objective function values are comparable for all four procedures.
Figure~\ref{fig:evo_proc_TCSD_DECREASING_order_constraints}
    shows the corresponding data profiles with a 5\% tolerance.
    
The gain made by the \intr and \hier procedures over the 
\eb and \pb
approaches are even more significant when the cost is low.  
For example, at an evaluation cost of 4,000, \hier and \intr have solved around 80\% of the instances, around twice as much as the procedures \eb and \pb.
The \intr procedure dominates all others on the totality of the graph,
    except that the \hier procedure is slightly above at the moment when a first feasible solution is produced.


\section{Discussion}
\label{sec-Discussion}
We have presented approaches that exploit 
    the flexibility of the \mads algorithm 
    to interrupt the evaluation sequence of the constraints as soon
    as one is shown to be infeasible.
Complex industrial problems, such as the problem PRIAD mentioned in introduction, will undoubtedly benefit from evaluation interruptions, and this without losing solution quality.
Numerical experiments suggest that the order in which the constraints are  evaluated affects the behaviour.
Future work include dynamically ordering the constraints, 
    as well as inserting the evaluation of the objective function
    in the sequence, and more importantly, applying this approach to real blackbox optimization problems.

\bibliographystyle{plain}
\bibliography{bibliography}

\begin{thebibliography}{10}

\bibitem{AuDe2006}
C.~Audet and J.E. {Dennis, Jr.}
\newblock {Mesh Adaptive Direct Search Algorithms for Constrained
  Optimization}.
\newblock {\em SIAM Journal on Optimization}, 17(1):188--217, 2006.

\bibitem{AuDe09a}
C.~Audet and J.E. {Dennis, Jr.}
\newblock {A Progressive Barrier for Derivative-Free Nonlinear Programming}.
\newblock {\em SIAM Journal on Optimization}, 20(1):445--472, 2009.

\bibitem{AuHa2017}
C.~Audet and W.~Hare.
\newblock {\em {Derivative-Free and Blackbox Optimization}}.
\newblock Springer Series in Operations Research and Financial Engineering.
  Springer, Cham, Switzerland, 2017.

\bibitem{AuLeDTr2018}
C.~Audet, S.~{Le~Digabel}, and C.~Tribes.
\newblock {The Mesh Adaptive Direct Search Algorithm for Granular and Discrete
  Variables}.
\newblock {\em SIAM Journal on Optimization}, 29(2):1164--1189, 2019.

\bibitem{belegundu1985study}
A.D. Belegundu and J.S. Arora.
\newblock {A study of mathematical programming methods for structural
  optimization. Part I: Theory}.
\newblock {\em International Journal for Numerical Methods in Engineering},
  21(9):1583--1599, 1985.

\bibitem{CoPa2019}
P.~C\^ot\'e and J.~Paquin.
\newblock {Optimization and Simulation Tools for Rio Tinto's Kemano Hydropower
  System}.
\newblock In {\em {Seventh International Megaprojects Workshop}}, Montr\'eal,
  June 2019.
\newblock Slides available at
  \url{https://sites.grenadine.uqam.ca/sites/megaprojectworskhop/fr/megaprojectworkshop/documents/get_document/39}.

\bibitem{FlLe02a}
R.~Fletcher and S.~Leyffer.
\newblock Nonlinear programming without a penalty function.
\newblock {\em Mathematical Programming}, Series A, 91:239--269, 2002.

\bibitem{garg2014solving}
H.~Garg.
\newblock {Solving structural engineering design optimization problems using an
  artificial bee colony algorithm}.
\newblock {\em Journal of Industrial and Management Optimization},
  10(3):777--794, 2014.

\bibitem{PRIAD_KoMeCoGaVoAlDeBl2021}
D.~Komljenovic, D.~Messaoudi, A.~C\^ot\'e, M.~Gaha, L.~Vouligny, S.~Alarie,
  A.~Dems, and O.~Blancke.
\newblock {Asset Management in Electrical Utilities in the Context of Business
  and Operational Complexity}.
\newblock In {\em {14th WCEAM Proceedings}}, 2021.

\bibitem{Le09b}
S.~{Le~Digabel}.
\newblock {Algorithm 909: NOMAD: Nonlinear Optimization with the MADS
  algorithm}.
\newblock {\em {ACM} Transactions on Mathematical Software}, 37(4):44:1--44:15,
  2011.

\bibitem{McCoBe79a}
M.D. McKay, R.J. Beckman, and W.J. Conover.
\newblock A comparison of three methods for selecting values of input variables
  in the analysis of output from a computer code.
\newblock {\em Technometrics}, 21(2):239--245, 1979.

\bibitem{MoWi2009}
J.J. Mor\'e and S.M. Wild.
\newblock {Benchmarking Derivative-Free Optimization Algorithms}.
\newblock {\em SIAM Journal on Optimization}, 20(1):172--191, 2009.

\end{thebibliography}
\end{document}